\input amstex
\documentstyle{amsppt}
\magnification1200
\tolerance=10000
\def\n#1{\Bbb #1}
\def\p{\Bbb C_{\infty}}
\def\r{\Bbb R_{\infty}}

\def\Ext{\hbox{Ext}}

\def\Ker{\hbox{Ker }}

\def\wt{\hbox{wt }}

\def\d{\Delta}

\def\e11{E_{11}}

\def\vf{\varphi}

\def\vf{\varphi}

\def\vi{v_\infty}
\def\g{\goth }

\def\om{\omega}
\def\th{\theta}

\topmatter
\title
Non-injectivity of the lattice map for non-mixed Anderson t-motives, and a result towards its surjectivity
\endtitle
\author
A. Grishkov, D. Logachev\footnotemark \footnotetext{E-mails: shuragri{\@}gmail.com; logachev94{\@}gmail.com (corresponding author)\phantom{****************}}
\endauthor
\thanks Thanks: The first author was supported by FAPESP (grant 2018/23690-6), 
CNPq(grants 307593/2023-1 and 406932/2023-9), and in accordance
 with the state task of the IM SB RAS, project FWNF-2022-003. 
\endthanks
\NoRunningHeads
\address
First author: Instituto de Matem\'atica e estatistica, 
Universidade de S\~ao Paulo. Rua de Mat\~ao 1010, CEP 05508-090, S\~ao Paulo, Brasil, 
\medskip
Omsk State University n.a. F.M.Dostoevskii, Russia,
\medskip
and Sobolev Institute of Mathematics, Omsk, Russia.
\medskip
Second author: Departamento de Matem\'atica, Universidade Federal do Amazonas, Manaus, Brasil.
\endaddress
\abstract Let $M$ be an uniformizable Anderson t-motive and $L(M)$ its lattice. First, we prove by an explicit construction that for the non-mixed $M$ the lattice map $M\mapsto L(M)$ is not injective. Second, we show that some lattices which do not belong to the set $L(M)$ of pure $M$, are lattices of non-pure $M$. This is a result towards surjectivity of the lattice map. The t-motives used in the proofs are non-pure t-motives of dimension 2, rank 3. Finally, we start calculations in order to answer a question whether all these t-motives are uniformizable, or not. 
\endabstract
\keywords Anderson t-motives; Lattice map \endkeywords
\subjclass 11G09 \endsubjclass
\endtopmatter

\NoRunningHeads

\document
{\bf 0. Introduction.} We shall consider in this paper only Anderson t-motives having the nilpotent operator $N$ equal to 0. Let $M$ be an  uniformizable t-motive of rank $r$, dimension $n$. Let $L(M)\subset \p^n$ be its lattice, it is isomorphic to $\n F_q[\th]^r$. Is the lattice map $M\mapsto L(M)$ an injection; a surjection? Unlike the case of abelian varieties where there is a 1 -- 1 correspondence between abelian varieties and lattices having a Riemann form, in the functional field case we know much less. 
\medskip
The present paper contains two results on the lattice map. First, we prove that for non-mixed t-motives the lattice map is not injective. Second, we show that some lattices which are not images of the lattice map of duals of Drinfeld modules, are lattices of some non-pure t-motives. 
\medskip
The technique of both results is similar to the one of [GL17]. We consider the simplest non-trivial case of t-motives of rank 3, dimension 2, defined by explicit equations (2.3), (3.5). We solve equations defining the lattices of these t-motives. 
\medskip
Let us explain what was known on the lattice map earlier. First, for the case of Drinfeld modules, i.e. for $n=1$, we have a 1 -- 1 correspondence between Drinfeld modules (they are all uniformizable) and lattices in $\p$ ([D76]). 
\medskip
Further, there is a result (see [GL17]) on a local surjectivity of the lattice map near $M=\g C_2^{\oplus n}$, where $\g C_2$ is the rank 2 Carlitz module (see 2.1 for its definition).
\medskip 
The most important result is the injectivity\footnotemark \footnotetext{This is only a rough statement.} of the lattice map on the set of mixed uniformizable Anderson t-motives ([HJ20], Theorem 3.34b). Recall that a t-motive is called mixed if it has a filtration such that its quotients are pure, of increasing weight. See [HJ20], Definition 3.5b for an exact definition of mixedness, and [G96], Definition 5.5.2 for the definition of purity. 
\medskip
The first result of the present paper: the condition of mixedness is necessary. We give two different examples of non-injectivity. First example: we give an explicit construction of a one-parametric family of non-isomorphic non-mixed $M(a)$, where $a\in \p$ is a parameter, such that all lattices $L(M(a))$ coincide. Second example: let $L$ be a fixed lattice. We show that there exist a pure $M_1$ and a non-pure $M_2$ (hence $M_1$, $M_2$ are not isomorphic) such that $L=L(M_1)=L(M_2)$. 
\medskip
 Let us describe the second result. By the theory of duality (see [GL07]), we have some information for the case $n=r-1$. Namely, let $L$ be a lattice of rank $r$ in $\p^{r-1}$. If $L^*$ --- the dual of $L$ --- exists, then it is a lattice of rank $r$ in $\p$. Hence, there exists a Drinfeld module $M$ such that $L(M)=L^*$. The dual of $M$ (in the meaning of [T95], [GL07]), denoted by $M^*$, always exists, it is a pure t-motive of rank $r$, dimension $r-1$. By the theory of duality ([GL07]), $L=L(M^*)$. Moreover, we have: 
\medskip
{\bf Theorem 0.1.} ([GL07], Corollary 8.4.) All pure t-motives of rank $r$, dimension $r-1$ are uniformizable. The lattice map is a 1 -- 1 correspondence between the set of pure t-motives of rank $r$, dimension $r-1$, and the set of lattices of rank $r$ in $\p^{r-1}$ having duals.

There is a natural 

{\bf Question 0.2.} Let $L$ be a lattice of rank $r$ in $\p^{r-1}$ such that its dual $L^*$ does not exist. Is $L$ a lattice of a t-motive (necessarily non-pure)? 
\medskip
Our second result (Theorem 3.4) is that for $r=3$, for many such $L$ the answer is "yes". This gives us evidence that the lattice map is surjective. 
 \medskip
Finally, we state the following question: what is the simplest example of a non-uniformizable t-motive $M$ (i.e. what are the minimal values of $r, \ n)$ of $M$? See Question 4.3 for more details. 
\medskip
The paper is organized as follows. In Section 1 we give necessary definitions and notations. Section 2 contains the proof of non-injectivity of the lattice map. Section 3 contains the proof of the second result. Section 4 contains some related  possibilities of further research. In Section 5 we start calculations to answer Question 4.3.
\medskip
{\bf 1. Definitions and notations.} We use standard definitions for Anderson t-motives. Let $q$ be a power of a prime, $\n F_q$ the finite field of order $q$, $\th$ a transcendental element, $\r:=\n F_q((1/\th))$ the finite characteristic analog of $\n R$. It has a valuation $v_\infty$ defined by $v_\infty(\th)=-1$. Let $\overline{\r}$ be an algebraic closure of $\r$. Let $\p:=\widehat{\overline{\r}}$ be the completion of $\overline{\r}$ with respect to (the only continuation of) $v_\infty$. It is the finite characteristic analog of $\n C$. 
\medskip
The Anderson ring $\p[T,\tau]$ is the ring of non-commutative polynomials over $\p$ in two variables $T$, $\tau$ with
the following relations (here $a\in \p$):
$$aT=Ta; \ \ \tau T=T\tau; \ \ \tau a=a^q\tau$$
It has subrings $\p[T]$, $\p\{\tau\}$.
\medskip
Let $A$ be a matrix with entries in $\p[T,\tau]$. We denote by $A^{(k)}$, where $k \in \n Z$, the matrix obtained by elevation of all coefficients of all entries of $A$ to the $q^k$-th power ($T$ and $\tau$ are not elevated to a power). 
\medskip
{\bf Definition 1.1.} An Anderson t-motive $M$ is a left $\p[T,\tau]$-module satisfying conditions:
\medskip
1.1.1. $M$ as a $\p[T]$-module is free of finite dimension (denoted by  $r$);
\medskip
1.1.2. $M$ as a $\p\{\tau\}$-module is free of finite dimension (denoted by  $n$);
\medskip
1.1.3. The action of $T-\th$ on $M/\tau M$ is nilpotent.
\medskip
Let $e_*:=\left(\matrix e_1\\ \dots \\ e_n  \endmatrix\right)$ be a basis of $M$ over $\p\{\tau\}$. In order to define $M$, we need to define the product $Te_*$:
$$Te_*=A_0e_*+A_1\tau e_*+...+A_k\tau^ke_*\eqno{(1.2)}$$
where $A_i\in M_{n\times n}(\p)$. Condition 1.1.3 is equivalent to $A_0=\th I_n +N$ where $N\in M_{n\times n}(\p)$ is nilpotent. We shall consider only $M$ having $N=0$. 
\medskip
Let $V=\p^n$ and $L\subset V$ is isomorphic to $\n F_q[\th]^r$. 
\medskip
{\bf Definition 1.3.} $L$ is called a lattice of rank $r$ if:
\medskip
{\bf 1.3.1.} The $\p$-span of $L$ is $V$;
\medskip
{\bf 1.3.2.} The $\r$-span of $L$ has dimension $r$ over $\n R_\infty$ (i.e. elements of a basis of $L$ over $\n F_q[\theta]$ are linearly independent over $\n R_\infty$).
\medskip
Two lattices $L_1\subset V_1$, $L_2\subset V_2$ are isomorphic if there exists a $\p$-linear isomorphism $\vf: V_1\to V_2$ such that $\vf(L_1)=L_2$. 
\medskip
{\bf 1.4. Lattice of a t-motive.} Let us fix a basis $e_*$ of $M$ over $\p\{\tau\}$, and let $A_i$ be from (1.2). The basis $e_*$ defines the following action of $T$ on $\p^n$:
$$T(Z)=\th Z+A_1Z^{(1)}+...+A_kZ^{(k)}$$
where $Z\in \p^n$ is a matrix column (recall that $N=0$, hence in our case $ A_0=\th I_n$). 
\medskip
The below theorems are proved in [D76] for the case $n=1$, and in [A86] for the case of any $n$.
\medskip
{\bf Theorem 1.4.1.} For a fixed $e_*$ there exists the only map $exp=exp_M: \p^n\to \p^n$ defined by the formula 
$$exp(Z)=Z+C_1Z^{(1)}+C_2Z^{(2)}+...\eqno{(1.4.1a)}$$
where $Z\in \p^n$ is a matrix column and $C_i\in M_{n\times n}(\p)$, making the following diagram commutative:
$$\matrix \p^n &\overset{exp}\to{\to} &\p^n &\\ \\ \theta\downarrow &&\downarrow & Z\mapsto T(Z) \\  \\ \p^n &\overset{exp}\to{\to} & \p^n &\endmatrix \eqno{(1.4.1b)}$$

{\bf Theorem 1.4.2.} $L(M):=\Ker exp$ is a $\n F_q[\th]$-submodule of $\p^n$ of dimension $\le r$.
\medskip
{\bf Definition 1.4.3.} $M$ is called uniformizable if the dimension of $L(M)$ as a $\n F_q[\th]$-module is $r$. 
\medskip
{\bf Theorem 1.4.4.} If $M$ is uniformizable then $L(M)$ is a lattice in $\p^n$. It is well-defined, i.e. if we change the basis $e_*$ by another basis $e'_*$ then we get an isomorphic lattice. 
\medskip
We shall consider $M$ of rank 3, dimension 2. Let us give explicit forms of (1.2) for them. First, let $M$ be pure. Its equation (1.2) is
$$T\left(\matrix e_1\\ e_2\endmatrix\right)=\left(\matrix \th &0 \\ 0&\th \endmatrix\right)\left(\matrix e_1\\ e_2\endmatrix\right)+\left(\matrix 0&-a_1\\1&-a_2 \endmatrix\right)\tau \left(\matrix e_1\\ e_2\endmatrix\right)+\left(\matrix 0&1\\0&0 \endmatrix\right)\tau^2\left(\matrix e_1\\ e_2\endmatrix\right)\eqno{(1.5)}$$
We denote it by $M_p(a_1, a_2)$ (subscript "p" means pure).
\medskip
{\bf Theorem 1.6.} $M_p(a_1, a_2)$ is pure, it is the dual of a Drinfeld module of rank 3 defined by the equation 
$$Te_1=\th \ e_1 + a_1 \tau e_1 + a_2 \tau^2 e_1 + \tau^3 e_1\eqno{(1.6.1)}$$
([T95], Section 5; [GL07]; [GL20], (12.2.2)). If $a_1, \ a_2$ are fixed then there are only finitely many $a'_1, \ a'_2$ such that $M_p(a_1, a_2)=M_p(a'_1, a'_2)$ ([D76]). All pure $M$ of rank 3, dimension 2 have this equation ([GL07]).
\medskip
Let now $M$ be a t-motive such that its equation (1.2) is
$$T\left(\matrix e_1\\ e_2\endmatrix\right)=\left(\matrix \th &0 \\ 0&\th \endmatrix\right)\left(\matrix e_1\\ e_2\endmatrix\right)+A\tau \left(\matrix e_1\\ e_2\endmatrix\right)+\left(\matrix 1&0\\0&0 \endmatrix\right)\tau^2\left(\matrix e_1\\ e_2\endmatrix\right)\eqno{(1.7)}$$
where $$A=\left(\matrix a_{11}&a_{12}\\ a_{21} &1\endmatrix\right)\eqno{(1.8)}$$ $M$ is not pure. We denote it by $M_{np}(A)$. 
\medskip
{\bf Conjecture 1.9.} All non-pure $M$ of rank 3, dimension 2 have an equation of the form (1.7). If $A$ of type (1.8) is fixed then there are only finitely many $A'$ of type (1.8) such that $M_{np}(A)=M_{np}(A')$.
\medskip
{\bf Corollary 1.10.} The set of non-pure $M$ of rank 3, dimension 2 is 3-dimensional. 
\medskip
We shall work with reducible $M$ defined by (1.7), see equations (2.3), (3.5).
\medskip
We see that $L$ is a functor (a map on the level of sets) from the set of uniformizable t-motives to the set of lattices, both up to isomorphisms. Is it injective or surjective? 
\medskip
Let us prove the result of non-injectivity: we construct explicitly a set of non-isomorphic uniformizable t-motives such that their lattices are isomorphic. 
\medskip
{\bf 2. Construction.} 
\medskip
We need more notations and definitions. First, the Carlitz module $\g C$ is a t-motive having $n=r=1$, it is defined by the formula $$Te=\th e+\tau e$$ where $e$ is the only element of a basis of $\g C$ over $\p\{\tau\}$. Let us denote $\th_{ij}:=\th^{q^i}-\th^{q^j}$ and for $j\ge 1$ $$c_j:=\frac1{\th_{j,j-1}\th_{j,j-2}\cdot...\cdot \th_{j1}\th_{j0}}$$
We have (Carlitz):
$$exp_\g C(z)=z+c_1z^q+c_2z^{q^2}+...$$
We denote by $\pi_1$ a generator of the lattice $L(\g C)$, it is unique up to multiplication by $\n F_q^*$. It is the finite characteristic analog of $2\pi i\in \n C$ which is defined up to $\pm1$.
\medskip
{\bf 2.1.} Further, the Carlitz module of rank 2 (denoted by $\g C_2$) is a t-motive having $n=1,\ r=2$ (i.e. it is a Drinfeld module of rank 2) defined by the formula $$Te=\th e+\tau^2 e$$
For even $j\ge2$ we denote $$c_{2,j}:=\frac1{\th_{j,j-2}\th_{j,j-4}\cdot...\cdot \th_{j2}\th_{j0}}$$
We have:
$$exp_{\g C_2}(z)=z+c_{2,2}z^{q^2}+c_{2,4}z^{q^4}+...$$
We denote by $\pi_2$ a generator of the lattice $L(\g C_2)$, it is unique up to multiplication by $\n F_{q^2}^*$. 

We choose and fix $\om\in \n F_{q^2}-\n F_{q}$, i.e. $(1,\om)$ is a basis of $\n F_{q^2}$ over $\n F_{q}$. We choose and fix $\pi_1$, $\pi_2$, hence $(\pi_2, \om\pi_2)$ is a basis of $L(\g C_2)$ over $\n F_{q}[\th]$. 
\medskip
Let $M_0:=\g C_2\oplus \g C$, and let $L_0:=L(M_0)\subset \p^2$ be its lattice. We have: elements
$$l_1:=\left(\matrix 0\\ \pi_1 \endmatrix\right), \ \ \ l_2:=\left(\matrix \pi_2\\ 0 \endmatrix\right), \ \ \ l_3:=\left(\matrix \om\pi_2\\ 0 \endmatrix\right)\eqno{(2.2)}$$ is a basis of $L_0$ over $\n F_{q}[\th]$.

Let $a\in \p$, and let $A$ from (1.8) be $\left(\matrix 0&a\\0&1 \endmatrix\right)$. We denote the corresponding $M(A)$ by $M(a)$. Its equation (1.7) is
$$T\left(\matrix e_1\\ e_2\endmatrix\right)=\left(\matrix \th &0 \\ 0&\th \endmatrix\right)\left(\matrix e_1\\ e_2\endmatrix\right)+\left(\matrix 0&a\\0&1 \endmatrix\right)\tau \left(\matrix e_1\\ e_2\endmatrix\right)+\left(\matrix 1&0\\0&0 \endmatrix\right)\tau^2\left(\matrix e_1\\ e_2\endmatrix\right)\eqno{(2.3)}$$
$M(a)$ is reducible: it enters in an exact sequence 
$$0\to \g C\to M(a)\to \g C_2\to0\eqno{(2.3.1)}$$
Arguments similar to [HJ20], Example 3.9 show that $M(a)$ is not mixed (if $a\ne0$), because $\wt(\g C)=1>\wt(\g C_2)=\frac12$. 
\medskip
We shall prove: 
\medskip
{\bf Proposition 2.4.} There exists a neighborhood $U$ of $0\in \p$ such that $\forall \ a\in U $  we have: $M(a)$ is uniformizable and $L(M(a))=L_0$.
\medskip
{\bf Proposition 2.5.} Conjecture 1.9 holds for these $M(a)$, i.e. for a fixed $a$ there are only finitely many $a'$ such that $M(a)=M(a')$. 
\medskip
This implies that the lattice map is not injective. 
\medskip
{\bf Proof of Proposition 2.4.} We need the notion of a Siegel matrix of a lattice. It is defined exactly like in the case of abelian varieties. Let $l_1, \dots, l_r$ be a basis of $L$ over $\n F_q[\th]$ such that $l_1, \dots, l_n$ is a basis of $V$ over $\p$. The Siegel matrix $S\in M_{(r-n)\times n}(\p)$ of $L$ with respect to a basis $l_1, \dots, l_r$ is defined by the formula $$\left(\matrix l_{n+1} \\ \dots \\ l_r  \endmatrix\right)=S \left(\matrix l_{1} \\ \dots \\ l_n  \endmatrix\right)$$
Two lattices $L_1$, $L_2$ are isomorphic iff there exist bases $l_{1*}, \ l_{2*}$ of $L_1$, $L_2$ such that their Siegel matrices $S_1$, $S_2$ are equal. 
\medskip
(2.2) implies that $S$ of $L_0$ is $(0, \om)$. Let us show that $\forall \ a\in\p$ a Siegel matrix of $L(M(a))$ is the same. $C_i$ from (1.4.1a) for $M(A)$ satisfy the recurrence relation
$$C_m=\frac{A_1C^{(1)}_{m-1}+A_2C^{(2)}_{m-2}}{\th_{m0}}$$ where $C_m=0$ for $m<0$, $C_0=I_2$, $A_1, \ A_2$ are from (1.7), i.e. for our case $A_1=\left(\matrix 0&a\\0&1 \endmatrix\right)$, $A_2=\left(\matrix 1&0\\0&0 \endmatrix\right)$. We get by induction that $$C_m=C_m(a)=\left(\matrix c_{2,m}&*\\0&c_m \endmatrix\right)\hbox{ for $m$ even, } C_m=C_m(a)=\left(\matrix 0&*\\0&c_m \endmatrix\right)\hbox{ for $m$ odd. } \eqno{(2.4.1)}$$ (see (3.7.1) for the formulas for (*), here we do not need them). 

For sufficiently small $a$ we have: $M(a)$ is uniformizable.\footnotemark \footnotetext{Most likely $M(a)$ is uniformizable for all $a$. We do not need this fact. See Remark 3.10a on uniformizability of $M_t(a)$, and Question 4.3 for a general problem.} We see that there is a basis $l_1(a), \ l_2(a), \ l_3(a)$ of $L(M(a))$ over $\n F_q[\th]$ such that $$ l_2(a)=l_2=\left(\matrix \pi_2\\ 0 \endmatrix\right), \ \ \ l_3(a)=l_3=\left(\matrix \om\pi_2\\ 0 \endmatrix\right), \hbox{ while }l_1(a)\ne l_1=\left(\matrix 0\\ \pi_1 \endmatrix\right).$$ We get that the Siegel matrix of $L(M(a))$ with respect to $l_1(a), \ l_2(a), \ l_3(a)$ is the same $(0,\om)$. $\square$
\medskip
{\bf Proof of Proposition 2.5.} We need a few facts on $GL_2(\p\{\tau\})$. Let 

$$X:=X_0+X_1\tau+... +X_k\tau^k\in M_{2\times2}(\p\{\tau\}),\eqno{(2.7)}$$ where $X_i= \left(\matrix x_{i11} & x_{i12} \\ x_{i21} & x_{i22}  \endmatrix\right)\in M_{2\times2}(\p).$ We have:
\medskip
$\exists \ Y\in M_{2\times2}(\p\{\tau\})$ such that $XY=1 \implies |X_0|\ne0 \implies X $ is not a zero divisor from both right and left $\implies YX=1$, i.e. $X\in GL_2(\p\{\tau\})$. 
\medskip
Further, if for $X\in GL_2(\p\{\tau\})$ we have $\forall \ i \ \ x_{i21}=0$ then $$\sum_{i=0}^k  x_{i11}\tau^i,  \  \sum_{i=0}^k  x_{i22}\tau^i\in \p\{\tau\}^*\implies \forall \ i >0 \ \ x_{i11}=x_{i22}=0.\eqno{(2.8)}$$

Let us assume that $M(a)=M(a')$. Multiplication by $T$ in $M(a')$ is given by 
$$T\left(\matrix e'_1\\ e'_2\endmatrix\right)=\left(\matrix \th &0 \\ 0&\th \endmatrix\right)\left(\matrix e'_1\\ e'_2\endmatrix\right)+\left(\matrix 0&a'\\0&1 \endmatrix\right)\tau \left(\matrix e'_1\\ e'_2\endmatrix\right)+\left(\matrix 1&0\\0&0 \endmatrix\right)\tau^2\left(\matrix e'_1\\ e'_2\endmatrix\right)\eqno{(2.9)}$$
An isomorphism between $M(a)$ and $M(a')$ is given by a matrix $X$ from (2.7) of change of basis (where $x_{***}$ are indeterminate coefficients), i.e. $$\left(\matrix e_1\\ e_2\endmatrix\right)=X\left(\matrix e'_1\\ e'_2\endmatrix\right).\eqno{(2.10)}$$ We substitute (2.10) to (2.3), and we multiply (2.9) by $X$ from the left. We get matrix equalities:
$$TX\left(\matrix e'_1\\ e'_2\endmatrix\right)=\left(\matrix \th &0 \\ 0&\th \endmatrix\right)X\left(\matrix e'_1\\ e'_2\endmatrix\right)+AX^{(1)}\tau \left(\matrix e'_1\\ e'_2\endmatrix\right)+\left(\matrix 1&0\\0&0 \endmatrix\right)X^{(2)}\tau^2\left(\matrix e'_1\\ e'_2\endmatrix\right)\eqno{(2.11)}$$
$$XT\left(\matrix e'_1\\ e'_2\endmatrix\right)=X\left(\matrix \th &0 \\ 0&\th \endmatrix\right)\left(\matrix e'_1\\ e'_2\endmatrix\right)+XA'\tau \left(\matrix e'_1\\ e'_2\endmatrix\right)+X\left(\matrix 1&0\\0&0 \endmatrix\right)\tau^2\left(\matrix e'_1\\ e'_2\endmatrix\right)\eqno{(2.12)}$$
where $A=\left(\matrix 0&a\\0&1 \endmatrix\right)$, $A'=\left(\matrix 0&a'\\0&1 \endmatrix\right)$.
Equality of coefficients at $\tau^m$ of (2.11) and (2.12) is
$$\th X_m+AX_{m-1}^{(1)}+\left(\matrix 1&0\\0&0 \endmatrix\right)X_{m-2}^{(2)}=\th^{q^m} X_m+X_{m-1}{A'}^{(m-1)}+X_{m-2}\left(\matrix 1&0\\0&0 \endmatrix\right)\eqno{(2.13)}$$
where $m=0, \dots, k+2$ and $X_i=0$ if $i\not\in 0,\dots,k$. 
\medskip
(2.13) shows that if $x_{m21}=x_{m-1,21}=0$ then $x_{m-2,21}=0$. By induction from up to down we get that $\forall \ m \  x_{m21}=0$. Hence, (2.8) implies that if $k\ge1$ then $X_k=\left(\matrix 0&*\\0&0 \endmatrix\right)$. 
Further, for $m=k+2$ (2.13) becomes 
$$X_k\left(\matrix 1&0\\0&0 \endmatrix\right)=\left(\matrix 1&0\\0&0 \endmatrix\right)X_k^{(2)}\iff x_{k11}\in \n F_{q^2}, \ x_{k12}=x_{k21}=0$$ 
This means that $k=0$, $X=X_0=\left(\matrix x_{011}&0\\0&x_{022} \endmatrix\right)$. Finally, (2.13) for $m=1$ is 
$$\left(\matrix 0&a\\ 0&1\endmatrix\right)X_0^{(1)}=X_0\left(\matrix 0&a'\\ 0&1\endmatrix\right), \hbox{ i.e. } \left(\matrix 0&x_{022}^q\ a\\ 0&x_{022}^q\endmatrix\right)=\left(\matrix 0&x_{011}\ a'\\ 0&x_{022}\endmatrix\right)  $$
This means that $x_{022}\in \n F_q$ and $M(a)$ is isomorphic to $M(a')$ iff $a'/a\in \n F_{q^2}$. $\square$
\medskip
{\bf 3. Lattices belonging to the image of the lattice map.} 
\medskip
We need a definition of the dual lattice. An invariant (i.e. do not depending of coordinates) definition is given for example in [GL07], Section 2; [GL20], (12.9.2); and [HJ20]. We shall need only the definition in terms of Siegel matrices. Namely, let $L\subset \p^n$ be a lattice, $l_*=(l_1,\dots, l_r)$ a $\n F_q[\th]$-basis of $L$ and $S$ the Siegel matrix of $L$ with respect to the basis $l_*$. 
\medskip
By definition, the dual lattice $L^*$ is a lattice of rank $r$ in $\p^{r-n}$ having a Siegel matrix $S^t$ in some basis. This notion is well-defined, i.e. it does not depend on a choice of $l_*$ and hence on a Siegel matrix. 
\medskip
Not all lattices have dual, because the condition (1.3.2) for $S^t$ is not always satisfied. Really, let $n=1$. A matrix $S=\left(\matrix s_1\\ \dots \\ s_{r-1}\endmatrix\right)\in M_{(r-1)\times 1}(\p)$ is a Siegel matrix of a lattice of rank $r$ in $\p$ iff 
$$1, s_1, \dots , s_{r-1}\hbox{ are linearly independent over } \r\eqno{(3.1)}$$ while its transposed $S^t=\left(\matrix s_1& \dots & s_{r-1}\endmatrix\right)\in M_{1\times (r-1)}(\p)$ is a Siegel matrix of a lattice of rank $r$ in $\p^{r-1}$ iff 
$$\exists \ i \hbox{ such that } s_i\not\in\r\eqno{(3.2)}$$
For $r>2$ (3.2) is weaker than (3.1), i.e. all lattices of rank $r$ in $\p$ have duals, but not all lattices of rank $r$ in $\p^{r-1}$ have duals. 
\medskip
{\bf Remark 3.3.} The same phenomenon occurs for Siegel matrices of other sizes. For example, let $S=\{s_{ij}\}\in M_{2\times 2}(\p)$ be a matrix such that $\vi(a_{11})\not\in \n Z$,  $\vi(a_{21})\not\in \n Z$,  $\vi(a_{11})-\vi(a_{21})\not\in \n Z$. Then obviously $S$ is a Siegel matrix of a lattice $(\n F_q[\th])^4\subset\p^2$. But if $1, \ a_{11}, \ a_{12}$ are linearly dependent over $\r$ then $S^t$ is not a Siegel matrix of a lattice. We think that it is necessary to modify the Definition 1.3 of a lattice, in order to get the surjectivity of the lattice map. 
\medskip
According Theorem 0.1, all lattices of rank $r$ in $\p^{r-1}$ having dual are lattices of pure t-motives. 
Let us answer Question 0.2. We consider the case $r=3$. We denote by $L(s_{11})$ the lattice having a Siegel matrix $S(s_{11}):=(s_{11}, \om)\in M_{1\times 2}(\p)$ (recall that $\om\in \n F_{q^2}-\n F_q$ is fixed). It is really a lattice, because $\om\not\in\r$, i.e. (3.2) holds. 
\medskip
We denote the field $\n F_{q^2}((\th^{-1})) \ = $ \{the $\r$-linear envelope of 1 and $\om$\} by $\n R_{\infty,2}$. 
If $s_{11}\in \n R_{\infty,2}$ then $L(s_{11})$ has no dual. 
\medskip
{\bf Theorem 3.4.} For all sufficiently small $s_{11}\in \p$ there exists a t-motive $M$ defined by (1.7) such that $L(M)=L(s_{11})$. 
\medskip
{\bf Proof.} Let $a\in \p$, and let $A$ from (1.8) be $\left(\matrix 0&0\\a&1 \endmatrix\right)$. We denote the corresponding $M(A)$ by $M_t(a)$ (the subscript "t" means "transposed"). Its equation (1.7) is
$$T\left(\matrix e_1\\ e_2\endmatrix\right)=\left(\matrix \th &0 \\ 0&\th \endmatrix\right)\left(\matrix e_1\\ e_2\endmatrix\right)+\left(\matrix 0&0\\a&1 \endmatrix\right)\tau \left(\matrix e_1\\ e_2\endmatrix\right)+\left(\matrix 1&0\\0&0 \endmatrix\right)\tau^2\left(\matrix e_1\\ e_2\endmatrix\right)\eqno{(3.5)}$$
These $M_t(a)$ enter in an exact sequence 
$$0\to \g C_2\to M_t(a)\to \g C\to0\eqno{(3.6)}$$ hence they are mixed. 
Their $C_m$ from (1.4.1a) are denoted by $C_{t,m}=C_{t,m}(a)$, they are transposed of $C_m(a)$ of (2.4.1):
$$C_{t,m}(a)=\left(\matrix c_{2,m}&0\\d_m(a)&c_m \endmatrix\right)\hbox{ for $m$ even, } C_{t,m}(a)=\left(\matrix 0&0\\d_m(a)&c_m \endmatrix\right)\hbox{ for $m$ odd. } \eqno{(3.7)}$$
The expression for $d_m(a)$ is the following, it can be easily found by induction: 

$$d_m(a)=\frac{[d_{m-1}(a)]^q}{\th_{m0}} \  \ \hbox { for $m$ even, } d_m(a)=\frac{a\cdot c_{2,m-1}^q}{\th_{m0}}+\frac{[d_{m-1}(a)]^q}{\th_{m0}} \  \ \hbox { for $m$ odd }$$

$$d_0(a)=0; \ d_1(a)=\frac{a}{\th_{10}};$$
$$d_2(a)=\frac{a^q}{\th_{21}\th_{20}};$$
$$d_3(a)=\frac{a}{\th_{31}\th_{30}}+\frac{a^{q^2}}{\th_{32}\th_{31}\th_{30}};\eqno{(3.7.1)}$$
$$d_4(a)=\frac{a^q}{\th_{42}\th_{41}\th_{40}}+\frac{a^{q^3}}{\th_{43}\th_{42}\th_{41}\th_{40}};$$
$$d_5(a)=\frac{a}{\th_{53}\th_{51}\th_{50}}+\frac{a^{q^2}}{\th_{53}\th_{52}\th_{51}\th_{50}} +\frac{a^{q^4}}{\th_{54}\th_{53}\th_{52}\th_{51}\th_{50}};$$
etc., we do not need its exact form. It is easy to find elements of a basis of $L(M_t(a))$, i.e. $Z=\left(\matrix z_1\\ z_2\endmatrix\right)\in \p^2$ which are solutions of 
$$\sum _{m=0}^\infty C_{t,m}(a)Z^{(m)}=0\eqno{(3.8)}$$
The solutions are $$l_1:=\left(\matrix 0\\ \pi_1\endmatrix\right),\ \ \  l_2:=\left(\matrix \pi_2\\ z_{22}\endmatrix\right),  \ \ \  l_3:=\left(\matrix \om\pi_2\\ z_{32}\endmatrix\right)\eqno{(3.9)}$$ To find $z_{22}$ and $z_{32}$ 
we denote $$D(a):=\sum _{m=1}^\infty d_m(a)\pi_2^{q^m}; \ \ \ D_\om(a):=\sum _{m=1}^\infty d_m(a)(\om\pi_2)^{q^m}$$
Both $D(a)$, $D_\om(a)$ are power series in $a$:
$$D(a)=\g d_0\ a\ +\ \g d_1\ a^q+\ \g d_2\ a^{q^2}+...$$
$$D_\om(a)=\bar \om\g d_0\ a\ +\ \om\g d_1\ a^q+\bar\om \g d_2\ a^{q^2}+...$$ where $\bar \om=\om^q$ is the conjugate of $\om$ and  $$\g d_0=\frac{\pi_2^q}{\th_{10}}+\frac{\pi_2^{q^3}}{\th_{31}\th_{30}}+\frac{\pi_2^{q^5}}{\th_{53}\th_{51}\th_{50}}+... , \  \g d_1=\frac{\pi_2^{q^2}}{\th_{21}\th_{20}}+\frac{\pi_2^{q^4}}{\th_{42}\th_{41}\th_{40}}+... , \g d_2=\frac{\pi_2^{q^3}}{\th_{32}\th_{31}\th_{30}}+... $$ 
Calculating $v_\infty(\g d_i)$ we see: 
\medskip
First, that $\g d_0\ne0$, because $v_\infty(\pi_2)=-q^2/(q^2-1)$ and $v_\infty(\g d_0)=-q/(q^2-1)$;
\medskip
Second, that $D(a)$ and $D_\om(a)$ converge for all $a\in \p$. 
\medskip
(3.8) becomes $$D(a)+exp_\g C(z_{22})=0; \ \ \ D_\om(a)+exp_\g C(z_{32})=0\eqno{(3.10)}$$

{\bf Remark 3.10a.} Since $exp_\g C$ is surjective, we get, as a by-product, that $M_t(a)$ is uniformizable $\forall \ a\in \p$. 
\medskip
Now and below we shall consider only sufficiently small $a$, such that the series $log_\g C(D(a))$, $log_\g C(D_\om(a))$ converge. We have: $$z_{22}=log_\g C(-D(a)), \ \ \ z_{32}=log_\g C(-D_\om(a))\eqno{(3.11)}$$
We denote the Siegel matrix of $M_t(a)$ corresponding to the basis (3.9), by $\g S(a)$. (3.9) and (3.11) imply that 
$$\g S(a)=S(\g s(a))=(\g s(a);\ \om\ )\eqno{(3.12)}$$
where $$\g s(a)=\frac{log_\g C(-D_\om(a))-\om\cdot log_\g C(-D(a))}{\pi_1}.$$
The function $a\mapsto \g s(a)$ is an additive power series in $a$ having a non-zero radius of convergence. It is non-zero: its first term is $$\frac{\g d_0(\om-\bar \om)}{\pi_1}a$$ and $\g d_0\ne0$. Hence, $a\mapsto \g s(a)$ is a local isomorphism in a neighborhood of 0. This completes the proof of Theorem 3.4: for $s_{11}$ near 0 we find $a$ such that $\g s(a)=s_{11}$. For this $a$ we have: $L(M_t(a))=L(s_{11})$. $\square$
\medskip
Let us consider two applications. 
\medskip
{\bf 3.13.} Let $s_{11}\in \n R_{\infty,2}$. The lattice $L(s_{11})$ has no dual. So, earlier we could not guarantee that it is the lattice of a t-motive (except the trivial case $s_{11}=0$). Theorem 3.4 tells us that for all sufficiently small $s_{11}$ there exists a t-motive such that $L(s_{11})$ is its lattice.
\medskip
{\bf 3.14.} Let a sufficiently small $s_{11}\in \p-\n R_{\infty,2}$. The lattice $L(s_{11})$ has dual. From one side, this means that $L(s_{11})$ is the lattice of the only one pure t-motive of rank 3, dimension 2 ( = dual of a Drinfeld module of rank 3). From another side, Theorem 3.4 tells us that $L(s_{11})$ is the lattice of a non-pure t-motive $M_t(a)$, as above. We get once again the result that the lattice map is not injective on the set of non-pure t-motives. 
\medskip
{\bf 4. Further questions.}
\medskip
{\bf Question 4.1.} Conjecture 1.9 implies that the set of $M(A)$ has dimension 3 over $\p$. Clearly the set of lattices of $M(A)$ has dimension 2 over $\p$, because it is defined by a $1\times 2$ Siegel matrix. Hence, we conjecture that the  fibers of the lattice map on the set of $M(A)$ have dimension 1. Conjecturally, they are analytic curves in the space $<a_{11}, a_{12}, a_{21}>$. Is it really true? What are these curves? 
\medskip 
{\bf Question 4.2.} Does exist a simply-described 2-dimensional subset $\g S_2$ of $A$ from (1.8) such that for any lattice $L\subset \p^2$ of rank 3 there exists $A\in \g S_2$ such that $L=L(M(A))$? I.e. this $\g S_2$ is a representative of the set of the above curves. 
\medskip
Existence of this $\g S_2$ will help us to prove the conjecture that the lattice map from the set of $M(A)$ defined by (1.7), to the set of lattices of rank 3 in $\p^2$, is surjective. 
\medskip
{\bf Question 4.3.} Are all $M(A)$ defined by (1.7) uniformizable? If not what are $h^1$, $h_1$ of these $M(A)$? 
\medskip
Recall that the simplest example of a non-uniformizable t-motive is given in [A86], Section 2.2 = [G96], Example 5.9.9. It has the minimal known values of $r=4$ and $n=2$. Whether there exist non-uniformizable t-motives having $r=3$, $n=2$ (they must be non-pure, i.e. conjecturally they are defined by the equation (1.7)), or not?

The technique developed in [GL21], [EGL22] will permit us to answer this question and to calculate $h^1$, $h_1$ of these $M(A)$ if they are $<3$ (recall that if $M$ is not uniformizable then it can happen that $h^1\ne h_1$). We start these calculations in Section 5. 
\medskip
$M(A)$ defined by (1.7) depend on 3 parameters while t-motives of rank 4, dimension 2 considered in [GL21], [EGL22] depend on 4 parameters. Hence, calculations for them will be simpler than the original calculations of [GL21], [EGL22]. We can start from computer calculations of $h^1$, $h_1$ of $M(A)$ for different $A$. 
\medskip
{\bf Question 4.4.} The set of $M(a)$ defined by (2.3) is $\Ext^1_{\p[T,\tau]}(\g C_2, \g C_1)$. We have: $\Ext^1_{\p[T,\tau]}(\g C_2, \g C_1)$ is a module over $Z(\p[T,\tau])=\n F_q[T]$. We showed in Proposition 2.5 that $\Ext^1_{\p[T,\tau]}(\g C_2, \g C_1)$ can be considered as having dimension 1 over $\p$. Hence, it is meaningful to consider a new invariant of a non-uniformizable t-motive $M$, namely the dimension of $\Ext^1_{\p[T,\tau]}(M, Z_1)$ (here $Z_1=\p\{T\}$ is from [G96], (5.9.22)) over $\p$. Recall that $M$ is uniformizable $ \iff \Ext^1_{\p[T,\tau]}(M, Z_1)=0$. As a $\n F_q[T]$-module, most likely, $\Ext^1_{\p[T,\tau]}(M, Z_1)$ is infinite-dimensional, but we can expect that over $\p$ the dimension is finite and gives us an invariant of $M$.
\medskip
It is possible to ask the same for Coker $\exp_M\subset \Ext^1_{\p[T,\tau]}(M, Z_1)$. 
\medskip
{\bf Remark 4.5.} It is not too difficult to find $v_\infty$ of the coefficients of the power series $\g s(a)$, and hence its Newton polygon. Hence, we can find explicitly the size of the set of $s_{11}$ that belong to the image of $\g s$. Maybe this set is even the whole $\p$. 
\medskip
But we want to prove the conjecture that the lattice map for $r=3, \ n=2$ is surjective. Therefore, we should not resprict ourselves by the Siegel matrices of the form $(s_{11}, \om)$, we should consider all matrices satisfying (3.2) and hence more general $M$ than the ones defined by (3.5). 
\medskip
{\bf 5. Lattices of t-motives defined by(1.7), (1.8).}
\medskip
Here we start to apply the methods of [GL21], [EGL22] to the non-pure t-motives $M$ of rank 3, dimension 2 defined by (1.7), (1.8). A basis of such $M$ over $\p[T]$ is 
$\{f_*\}=\left(\matrix f_1\\ f_2\\ f_3\endmatrix\right)=\left(\matrix e_1\\ e_2\\ \tau e_1\endmatrix\right)$. The action of $\tau$ on $\{f_*\}$ is defined by a matrix $Q\in \p[T]$: $$\tau f_*=Q f_* \hbox { where in our case } Q=\left(\matrix 0&0&1\\ 0&T-\th&-a_{21}\\ T-\th & -a_{12}(T-\th)&-d \endmatrix\right)$$ where $d=|A|=a_{11}-a_{12}a_{21}$. In order to find $L(M)$, we must solve a system $$QX=X^{(1)}\eqno{(5.1)}$$ where $X=\left(\matrix X_1\\ X_2\\X_3\endmatrix\right)\in \p\{T\}^3$. Hence, (5.1) is 
$$X_3=X_1^{(1)}\eqno{(5.2)}$$
$$(T-\th)X_2-a_{21}X_3=X_2^{(1)}\eqno{(5.3)}$$
$$(T-\th)X_1-a_{12}(T-\th)X_2-dX_3=X_3^{(1)}\eqno{(5.4)}$$
Let us make eliminations. (5.2) --- (5.4) imply
$$(T-\th^q)X_1^{(1)}-a_{12}^q(T-\th^q)X_2^{(1)}-d^qX_1^{(2)}-X_1^{(3)}=0\eqno{(5.5)}$$
$$X_1^{(1)}=\frac{T-\th}{a_{21}}X_2-\frac1{a_{21}}X_2^{(1)}\eqno{(5.6)}$$
$$X_1^{(2)}=\frac{T-\th^q}{a_{21}^q}X_2^{(1)}-\frac1{a_{21}^q}X_2^{(2)}\eqno{(5.7)}$$
$$X_1^{(3)}=\frac{T-\th^{q^2}}{a_{21}^{q^2}}X_2^{(2)}-\frac1{a_{21}^{q^2}}X_2^{(3)}\eqno{(5.8)}$$
Substituting (5.6) --- (5.8) to (5.5) we get
\medskip
$$\frac1{a_{21}^{q^2}}X_2^{(3)}+(-\frac{T}{a_{21}^{q^2}}+\frac{d^q}{a_{21}^q}+\frac{\th^{q^2}}{a_{21}^{q^2}})X_2^{(2)}-(T-\th^q)(\frac1{a_{21}}+a_{12}^q+\frac{d^q}{a_{21}^q})X_2^{(1)}+$$ $$+\frac{(T-\th^q)(T-\th)}{a_{21}}X_2=0\eqno{(5.9)}$$

We can introduce new variables $$u=\frac1{a_{21}}+a_{12}^q+\frac{d^q}{a_{21}^q}=\frac1{a_{21}}+\frac{a_{11}^q}{a_{21}^q}$$ $$v=\frac{d^q}{a_{21}^q}+\frac{\th^{q^2}}{a_{21}^{q^2}}.$$ If $a_{21}\ne 0$ then the set of $(a_{21}, a_{11}, a_{12})$ is in 1 -- 1 correspondence with the set of $(a_{21}, u, v)$. (5.9) has the form 
$$ \frac1{a_{21}^{q^2}}X_2^{(3)}+(-\frac{T}{a_{21}^{q^2}}+v)X_2^{(2)}-(T-\th^q)\ u\ X_2^{(1)}+\frac{(T-\th^q)(T-\th)}{a_{21}}X_2=0\eqno{(5.10)}$$
Continuing the calculations similar to the ones of [GL21], [EGL22], we shall find $h_1$, $h^1$ of these $M$.
\medskip
{\bf References}
\medskip
[A86] Anderson, G.W. t-motives. Duke Math. J., 1986, vol. 53, No. 2, p. 457 -- 502. 
\medskip
[D76] V.G. Drinfeld, Elliptic modules. Math. USSR Sb. 4 (1976) 561 -- 592.
\medskip
[EGL22] S. Ehbauer, A. Grishkov, D. Logachev. Calculation of $h^1$ of some Anderson t-motives. Journal of Algebra and its applications, 2022, vol. 21, no. 1, Paper No. 2250017. 
https://arxiv.org/pdf/2006.00316.pdf 
\medskip
[G96] Goss, D. Basic structures of function field arithmetic. Springer-Verlag, Berlin, 1996. xiv+422 pp.
\medskip
[GL07] Grishkov A., Logachev, D. Duality of Anderson t-motives. 2007. 

http://arxiv.org/pdf/math/0711.1928.pdf
\medskip
[GL17] Grishkov A., Logachev, D. Lattice map for Anderson t-motives: first approach. J. of Number Theory. 2017, vol. 180, p. 373 -- 402.

http://arxiv.org/pdf/math/1109.0679.pdf
\medskip
[GL20] Grishkov A., Logachev, D. Introduction to Anderson t-motives: a survey. 2020. https://arxiv.org/pdf/2008.10657.pdf
\medskip
[GL21] Grishkov A., Logachev D., $h^1\ne h_1$ for Anderson t-motives. J. of Number Theory. 2021, vol. 225, p. 59 -- 89. 
https://arxiv.org/pdf/1807.08675.pdf 
\medskip
[HJ20] Hartl U.; Juschka A.-K. Pink's theory of Hodge structures and the Hodge conjecture over function fields. 
"t-motives: Hodge structures, transcendence and other motivic aspects", Editors G. B\"ockle, D. Goss, U. Hartl, M. Papanikolas, 
European Mathematical Society Congress Reports 2020, and https://arxiv.org/pdf/1607.01412.pdf
\medskip
[T95] Taguchi, Yuichiro. A duality for finite t-modules. J. Math. Sci. Univ. Tokyo
2 (1995), no. 3, 563 -- 588.
\medskip
\enddocument